\input amstex
\documentstyle{amsppt}
\input epsf.tex
%
%
\def\C{{\Bbb C}}
\def\R{{\Bbb R}}
\def\Z{{\Bbb Z}}

\def\N{{\Bbb N}}

%
%

\def\Tors{\mathop{\roman{Tors}}\nolimits}

\def\Im{\mathop{\roman{Im}}\nolimits}

\def\Int{\mathop{\roman{Int}}\nolimits}

\def\abs{\mathop{\roman{abs}}\nolimits}

\def\conj{\mathop{\roman{conj}}\nolimits}

\def\Spin{\mathop\roman{Spin}\nolimits}
\def\Spinc{\mathop\roman{Spin^{\C}}\nolimits}
\def\SO{\mathop\roman{SO}\nolimits}

\def\Ske{\mathop\roman{Ske}\nolimits}

%
%
\def\a{\alpha}

\def\s{\sigma}
\def\om{\omega}

%

%
%
%
\let\ge\geqslant
\let\le\leqslant
\def\la{\langle}
\def\ra{\rangle}

\def\Cp#1{\C\roman P^{#1}}

\def\Rp#1{\R\roman P^{#1}}
%
%
\let\tm\proclaim
\let\endtm\endproclaim

\let\rk=\remark
\let\endrk=\endremark

%
%
%
%

\def\SS{\Sigma}

\def\SC{\Spin^{\C}}
\def\SW{\mathop{SW}\nolimits}

\def\Wedge{\Lambda}

\def\A{\Bbb A}

\def\OS{\mathop{OS}\nolimits}
\def\SCn{\Spin^{\C}_n}

\def\ss{\frak s}
\def\ttt{\frak t}
\def\rr{\frak r}
\def\abs{\mathop{abs}\nolimits}

\def\mix{\roman{mix}}

\hsize=388pt \vsize=572pt
 \baselineskip 20pt \magnification\magstep1 \nologo \NoRunningHeads
\NoBlackBoxes

\topmatter
\title
Relative Seiberg-Witten and Ozsvath-Szabo 4-dimensional
invariants with respect to embedded surfaces
\endtitle
\author Sergey Finashin
\endauthor
\endtopmatter


\document
\heading \S1. Introduction \endheading

\subheading{1.1. The subject} In this paper we define and study
the relative Seiberg-Witten (SW) invariant and an analogous
relative Ozsv\'ath-Szab\'o (OS) invariant of pairs $(X,\SS)$,
where $\SS$ is a surface of genus $g>1$ embedded in a $4$-manifold
$X$.
This kind of invariant in case of $g=1$ was constructed by Taubes
\cite{T5}. It was observed in \cite{T5} and then ingeniously used
in \cite{FS1} and \cite{FS2} that the relative SW invariant of
Taubes for $(X,\SS)$ is practically reduced to the absolute SW
invariant of the fiber sum of $X$ with $E(1)$.

This property suggested an elementary definition of a similar
invariant for $(X,\SS)$, $g>1$, as the pull-back of the absolute
invariant of a certain fiber sum $X\#_\SS W$. The same approach
works for the OS invariants. Our goal is to show that the choice
of $W$ is not essential provided it admits a relatively minimal
Lefschetz fibration $W\to S^2$ with a fiber $\SS$ and has
$H_1(W)=0$. We obtain also a package of properties for these
relative invariants which is analogous to the package of
properties of the absolute invariants.

Briefly speaking, if a basic $\SC$ structure, $\ss$, in SW (or in
OS) theory is extremal with respect to the adjunction inequality
for $\SS$, that is belongs to the image,
$\SC(X,\ss_\SS)=\{x\in\SC(X)|c_1(\ss)=\chi(\SS)+\SS^2\}$, of the
forgetful map $\abs_{X,\SS}\:\SC(X,\SS)\to\SC(X)$, then
$\abs_{X,\SS}^{-1}(\ss)$ contains several relative basic
structures, $\rr_i\in\SC(X,\SS)$, and the SW (as well as OS)
invariant of $\ss$ splits into a sum of the relative invariant of
$\rr_i$. The relative invariants carry a package of properties
analogous to those of the absolute invariants.

 The manifolds that we consider in what follows
are smooth, oriented, connected, and closed, unless we state
otherwise. We suppose also that the surface $\SS\subset X$ has
genus $g>1$ and either essential with self-intersection $\SS^2=0$,
or has $\SS^2>0$ (in the latter case we blow it up to obtain
$\SS^2=0$). This implies in particular that $b_2^+(X)\ge1$.

 \subheading{1.2. The Seiberg-Witten and the Ozsv\'ath-Szab\'o invariants}
In its simplest version, the SW invariant of a $4$-manifold $X$ is
a function on the set of $\SC$ structures, $\SW_X\:\SC(X)\to\Z$,
which takes non-zero values only at a finite set of
$\ss\in\SC(X)$, called {\it the basic structures}, whose degree
$d(\ss)=\frac14(c_1^2(\ss)-(2\chi(X)+ 3 \s(X))$ is zero. The
corresponding Ozsv\'ath-Szab\'o invariant (extracted as the
reduced form of $\Phi_{X,\ss}$ in \cite{OS4},\S4) is an analogous
function $\OS_X\:\SC(X)\to\Z/\pm$ well-defined up to sign.
The sign of $\SW_X(\ss)$ depends on the choice of an orientation
of $H^1(X;\R)\oplus H^2_+(X;\R)$ called the {\it homology
orientation}, which is supposed to be fixed. In a bit special case
of $b_2^+(X)=1$ the invariants $\SW_X$, $\OS_X$ depend on some
additional data that must be fixed: $\SW_X$ depends on an
orientation of the line $H^2_+(X)$ (see \cite{KM}, \cite{MST}, or
\cite{T1}), while $\OS_X$, according to \cite{OS4, Prop. 2.6},
depend on the choice of an isotropic line in $H_2(X)$.

One can consider also a more refined and sophisticated version of
$\SW_X$ and $\OS_X$, which may take non-zero values for
$\ss\in\SC(X)$ of degree $d(\ss)>0$. These values are
homomorphisms $\SW_{X,\ss}\:\A_X\to \Z$, $\OS_{X,\ss}\:\A_X\to
\Z/\pm$, supported in the set of homogeneous elements of degree
$d(\ss)$ from the graded ring
$\A_{X}=\Wedge(H_1(X)/\Tors)\otimes\Z[U]$, where $\Wedge$ stands
for the exterior algebra and the grading is defined on the
generators, so that $U$ has degree $2$ and $\a\in
H^1(X)\smallsetminus\{0\}$ have degree $1$.
 The dual consideration, which is more convenient for us,
 interprets the refined version of $\SW_X$ as a map $\SC(X)\to
\A_{X}^*=\Wedge(H^1(X))\otimes\Z[U]\cong \Wedge(H^1(X;\Z[U]))$,
such that $\SW_X(\ss)$ is homogeneous of degree $d(\ss)$. Reducing
the values of $\OS_X$ modulo $2$, we obtain a similar map
$\SC(X)\to \A_{X}^*\otimes\Z/2$.

We will use notation $S_X\:\SC(X)\to R_X$ for  any of the
invariants $\SW_X$ or $\OS_X$, either in the refined or in the
reduced form. $R_X$ here is a ring $\A_{X}^*$ or
$\A_{X}^*\otimes\Z/2$ in case of the refined forms of $\SW_X$ or
$\OS_X$. In case of the reduced forms, $R_X$ is just $\Z$ or
$\Z/2$. The relative version, $R_{X,\SS}$, of this ring is
obtained by replacing $H^1(X)$ by $H^1(X,\SS)$ in the definition
(so that $R_{X,\SS}=R_X$ for the reduced forms of SW and OS
invariants).

The formal sum 
$$\underline{S}_{X}=\sum_{\ss\in\SC(X)} S_X(\ss)\centerdot\ss$$ 
can be considered as an element of the
principal (affine) module $R_{X}[\SC(X)]$ over the group ring
$R_{X}[H^2(X)]$.

\rk{Remark} All the constructions and the results obtained below
for SW and OS invariants concern in fact any function $S_X$
satisfying a few basic properties of SW and OS invariants, namely
A1--A5 formulated in \S3.
\endrk

\subheading{1.3. Definition of the relative invariant $S_{X,\SS}$}
Let $\SC(X,\SS)$ denote the set of relative $\SC$ structures and
$\abs_{X,\SS}\:\SC(X,\SS)\to\SC(X)$ the forgetful map.
Gluing of relative $\SC$ structures in a fiber sum $X_+ \#_\SS
X_-$ (see 2.11) yields
$$
\vee\:\SC(X_+,\SS)\times\SC(X_-,\SS)\to\SC(X_+\#_\SS X_-,\SS),\ \
(\rr_+,\rr_-)\mapsto \rr_+\vee\rr_-
$$
whose composition with the forgetful map gives
$$\#_\SS\:\SC(X_+,\SS)\times\SC(X_-,\SS)\to\SC(X_+\#_\SS X_-),\ \
(\rr_+,\rr_-)\mapsto \rr_+\#_\SS\rr_- .
$$

Choose any relatively minimal Lefschetz fibration $W\to S^2$ with
a fiber $\SS$ and $b_1(W)=0$ and denote by
$\rr_{W,\SS}\in\SC(W,\SS)$ its canonical relative $\SC$ structure
of the Lefschetz fibration introduced in 2.7. If $\SS^2=0$, then
we define for any $\rr\in\SC(X,\SS)$
$$S_{X,\SS}(\rr)=S_{X\#_\SS W}(\rr\#_\SS \rr_{W,\SS})$$

If $\SS^2>0$, then we blow up $X$ at points of $\SS$ to obtain
$\hat X$, with $\SS^2=0$, and let
$$S_{X,\SS}(\rr)=S_{\hat X,\SS}(\hat\rr)$$
where $\hat\rr$ is the image of $\rr$ under the natural map
$\SC(X,\SS)\to\SC(\hat X,\SS)$ (see 2.10).

We will let $$\underline{S}_{X,\SS}=\sum_{\rr\in\SC(X,\SS)}
S_{X,\SS}(\rr)\centerdot\rr\in R_{X,\SS}[\SC(X,\SS)]$$

\rk{Remarks}\roster\item Note that $b_2^+(X\#_\SS W)>1$, so
$S_{X\#_\SS W}$ is well defined. \item The differential type of
$X\#_\SS W$ may depend in principle on the framing of $\SS$ in $X$
and $W$, so proving that $S_{X,\SS}$ is independent of $W$ implies
also independence of the framings.
 \item In the case $S_X=\SW_X$ we
should take care of the homology orientation for $X\#_\SS W$. It
is determined by the given homology orientation of $X$ and the
canonical symplectic homology orientation of $W$ (defined in
\cite{T4}) following the rule described in \cite{MST} after a
modification, which is just such an alternation of the homology
orientation which eliminates the sign $(-1)^{b(M,N)}$ that appears
in the product formula of \cite{MST}. In the other words, with
such a homology orientation the product formula will look like A4,
in \S3 below.

Note that such a homology orientation in the fiber sums is
preserved by the natural diffeomorphisms $X\#_\SS Y\cong Y\#_\SS
X$ and $(X\#_\SS Y)\#_\SS Z\cong X\#_\SS(Y\#_\SS Z)$. In the case
of symplectic pairs $(X,\SS)$ and $(Y,\SS)$, the symplectic
homology orientations in $X$ and $Y$ induce the symplectic
homology orientation of $X\#_\SS Y$.

If $\SS^2>0$, then we choose the homology orientation of $\hat X$
induced by that of $X$.
\endroster
\endrk

\subheading{1.4. The properties of $S_{X,\SS}$}
 \tm{1.4.1. Theorem} The invariant
$S_{X,\SS}\:\SC(X,\SS)\to R_{X,\SS}$ is independent of the choice
of $W$ and has the following properties.
 \roster
 \item  Finiteness of the set
$B_{X,\SS}=\{\rr\in\SC(X,\SS)\,|\,S_{X,\SS}(\rr)\ne0\}$.
 \item The blow-up relation  $S_{\hat X,\SS}(\hat\rr)=S_{ X,\SS}(\rr)$
$($if $\SS^2>0$ in $X$ $)$.
 \item
The conjugation symmetry
 $S_{X,\SS}=\pm S_{X,-\SS}\circ\conj$,
where $-\SS$ is $\SS$ with the opposite orientation, and
$\conj\:\SC(X,\SS)\to\SC(X,-\SS)$ the conjugation involution.
 \item
Normalization: if $X$ admits a relatively minimal Lefschetz pencil
with a fiber $\SS$ and is endowed with the canonical homology
orientation of a symplectic manifold (for a symplectic structure
compatible with the pencil), then there is only one basic relative
$\SC$ structure, $\rr_{X,\SS}\in\SC(X,\SS)$ (the canonical
structure of the pencil, see 2.6) and $S_{X,\SS}(\rr_{X,\SS})=1$.
\item Splitting formula relating the absolute and the relative SW
invariants:
$$S_X(\ss)=\sum_{r\in\abs_{X,\SS}^{-1}(\ss)}S_{X,\SS}(\rr),$$
for any $\ss\in\SC(X)$ such that $c_1(\ss)[\SS]=\chi(\SS)+\SS^2$.
 \item
The product formula for a fiber sum $X=X_+\#_{\SS} X_-$ says:
$$\underline S_{X,\SS}=(\underline S_{X_-,\SS})(\underline S_{X_+,\SS})
$$
More explicitly, this means that for any $\rr\in\SC(X,\SS)$ 
$$S_{X,\SS}(\rr)=\sum_{\rr_+\vee\rr_-=\rr}S_{X_+,\SS}(\rr_+)
S_{X_-,\SS}(\rr_-)$$
 where $\rr_\pm$ are varying in $\SC(X_\pm,\SS)$.
 Equivalently, one can write it as
$$
S_{X,\SS}(\rr_+\vee\rr_-)=\sum_{k\in\Z}S_{X_+,\SS}(\rr_++k\s_+)
 S_{X_-,\SS}(\rr_--k\s_-)
$$
where $\s_\pm\in H^2(X_\pm,\SS)$ is dual to the fundamental class
of $\SS$ shifted inside $X_\pm\smallsetminus\SS$.
 \item
Adjunction inequality for $\rr\in B_{X,\SS}$ and a membrane
$F\subset X$, with the connected complement
$\SS\smallsetminus\partial F$ (note that positivity of $F^2$ is
not required)
$$ -\chi(F)\ge F^2+|\rr[F]|.$$
\endroster
\endtm

 The relations (4) and (5) of Theorem 1.4.1 imply together
\tm{Corollary 1.4.2} For any $\ss\in\SC(X)$ such that
$c_1(\ss)[\SS]=\chi(\SS)+\SS^2$
$$S_X(\ss)=\sum\Sb\rr_+\#_\SS\rr_-=\ss\endSb S_{X_+,\SS}(\rr_+)
S_{X_-,\SS}(\rr_-)$$

\endtm

\rk{Remarks} \roster\item The sign ``$\pm$'' in Theorem 1.4.1(3)
is $(-1)^{\frac14(\sigma(X)+\chi(X))}$, like for the absolute
invariant $S_X$. \item The product of $\underline{S}_{X_\pm,\SS}$ in
Theorem 1.4.1(6), is induced by the natural affine map
$$\SC(X_+,\SS)\times \SC(X_-,\SS)\cong \SC(X,M)\to \SC(X,\SS)$$
associated with the corresponding cohomology homomorphisms.
 \item
If $b_2^+(X)=1$, then the splitting formula 1.4.1(5) should be
applied to $\SW_X(\ss)$ defined by fixing the $[\SS]$-positive
orientation of the line $H_2^+(X)$, and $\OS_X(\ss)$ is defined by
fixing the line spanned by $[\SS]$ (because these are the choices
involved into the corresponding product formulae).
 \item Theorem
1.4.1 does not mention some more straightforward properties, which
do not involve $\SS$ , for example, the adjunction inequality for
a closed surface $F$, the blowup formula at a point $x\notin\SS$,
and the product formula with respect to an additional surface
$\SS'$ in $X_\pm$ disjoint from $\SS$. All these properties are
formulated exactly like in the case of the absolute invariants
$S_X$ and proved by giving an obvious reference to the case of
absolute invariants of the corresponding fiber sums.
 \item The invariant $S_{X,\SS}$ can be
defined similarly for a multi-component surface $\SS$, as one can
take fiber sums with auxiliary Lefschetz fibrations along all the
components of $\SS$. The properties of such invariants are
analogous to those formulated in Theorem 1.4.1, and the proofs
just repeat the arguments in \S4.
 \item The case of genus $g=1$ is a bit special, mainly
because the product formula in this case looks different.
Nevertheless, the same definition for $S_{X,\SS}$ can be given for
$g=1$, and all the properties except the splitting formula (5) in
Theorem 1.4.1 still hold. In the case of $S_X=\SW_X$, this follows
from the results of Taubes \cite{T5}, except for the property (7)
(not discussed in \cite{T5}), which is proved by the same
arguments as in the case $g>1$.
\endroster
\endrk

\subheading{1.5. Application: the genus estimate for membranes} By
definition, a membrane on a surface $\SS$ in $X$ is a compact
surface $F\subset X$ with the boundary $\partial F=F\cap\SS$, at
no point of which $F$ is tangent to $\SS$. The self-intersection
index $F^2$ is defined with respect to the normal framing along
$\partial F$ which is tangent to $\SS$. The number $\rr[F]$
(evaluation of $\rr\in\SC(X,\SS)$ on $F$) is defined in 2.6.
Throughout the paper we suppose that membranes are connected and
oriented, although the adjunction inequality holds as well for
disconnected membranes, which follows from additivity of
$\chi(F)$, $F^2$ and $\rr[F]$.

The adjunction inequality for membranes implies for instance the
minimal genus property for symplectic and Lagrangian membranes in
symplectic manifolds, namely

\tm{1.5.1. Corollary} Assume that $X$ is a symplectic
$4$-manifold, $\SS\subset X$ is an essential surface with
$\SS^2\ge0$, $g(\SS)>1$, and $F\subset X$ is a membrane on $\SS$
with the connected complement $\SS\smallsetminus L$ of the
boundary $L=\partial F$. Assume furthermore that either $\SS$ is
symplectic and $F$ is Lagrangian, or vice versa, $\SS$ is
Lagrangian and $F$ is symplectic. Then for any membrane,
$F'\subset X$, which has the same boundary $L=\partial F'$,
defines the same normal framing for $\SS$ along each component of
$L$, and realizes the same class $[F',\partial F']=[F,\partial
F]\in H_2(X,\SS)$ we have $g(F')\ge g(F)$.
\endtm

\demo{Proof} The adjunction inequality 1.4.1(7) becomes an
equality for such $X$, $\SS$,  $F$, and for the canonical $\SC$
structure $\rr\in\SC(X,\SS)$ (symplectic or Lagrangian, depending
on the case considered). The assumptions formulated for $F'$ imply
that $(F')^2=F^2$ and $\rr[F']=\rr[F]$, so 1.4.1(7) yields the
required estimate for $g(F')$. \qed\enddemo

Another example of application of 1.4.1(7) is orthogonality of the
relative basic classes, $\rr\in B_{X,\SS}$, to the $(-1)$-disc
membranes.

\tm{1.5.2. Corollary} If $D\subset X$ is a $(-1)$-disc membrane on
a surface $\SS\subset X$ as above, then $\rr[D]=0$ for all $\rr\in
B_{X,\SS}$. \qed\endtm

\subheading{1.6. On the calculation of the invariants $S_{X,\SS}$}
The first observation concerns vanishing of the relative
invariants $S_{X,\SS}$ if $\SS$ is not a minimal genus surface in
its homology class, because non-vanishing would
 contradict to the adjunction inequality applied in
$X\#_\SS W$ to the surface $\SS'$ homologous to $\SS$ but of a
smaller genus. This argument does not work if $\SS'$ cannot be made
disjoint from $\SS$ in $X$ 
(although the author does not know such examples, in which it really cannot).

Another example of calculation of $S_{X,\SS}$ is
contained in Theorem 1.4.1(4).
It is an interesting question if the normalization property 1.4.1(4) holds as well
in a more general setting, namely for symplectic
relatively minimal pairs $(X,\SS)$, such that $[\SS]\in H_2(X)$ is a primitive class.

One more important example of calculation
of $\SW_{X,\SS}$ can be extracted from the work \cite{FS2}, where
these invariants appeared to distinguish the embeddings of
surfaces obtained by the rim-surgery from $\SS$. In fact, the
results of \cite{FS2} mean that the invariant $\SW_{X,\SS}$ is
multiplied by the Alexander polynomial after performing a rim knot
surgery.
More precisely, assume that $\ell\subset\SS$ is a simple closed
curve, $K\subset S^3$ is a knot and $\SS_{K,\ell}\subset X$ is a
surface obtained from $\SS\subset X$ by rim surgery along $\ell$
using $K$ as a pattern. 

\tm{1.6.1. Theorem}
$\underline{\SW}_{X,\SS_{K,\ell}}=\Delta_K(\delta([\ell]^*))\underline{\SW}_{X,\SS}$.
\qed\endtm

Here $[\ell]^*\in H^1(\SS)$ is dual to $[\ell]\in H_1(\SS)$,
$\delta$ is the boundary map $H^1(\SS)\to H^2(X,\SS)$, $\Delta_K$
is the Alexander polynomial in the symmetrized form, and
$\Delta_K(\delta([\ell]^*))$ is considered as an element of the
group ring $\Z[H^2(X,\SS)]\subset R_{X,\SS}[H^2(X,\SS)]$.

\subheading{1.7. The structure of the paper} In \S2 we give a
brief summary on the absolute $\SC$ structures and develop some
calculus of the relative structures that is used in \S4.

In \S3, we recall some fundamental properties of SW and OS
invariants that are required to construct their relative versions.
Although mostly well-known, these properties appear in literature
in various settings, not always in the form convenient for us, so
we give some comments and references. In the core section, \S4, we
prove Theorem 1.4.1.

In \S5, we discuss some generalizations of the invariant
$S_{X,\SS}$. Its version, $S_{X,\SS,K}$, discussed in 5.1, depends
on a subgroup $K\subset H^1(\SS)$, which suggests an analogy with
the invariant in \cite{CW}. A generalization of $S_{X,\SS}$ in
5.2, which we presented for simplicity only in the case of OS
invariants, is a relative version of the invariant $F^{\mix}_{X,\ss}$
from \cite{OS3} (here instead of relativization with respect to a surface
$\SS\subset X$ we consider relativization with respect 
to a boundary component of $X$).

\subheading{1.8. Acknowledgements} The proof of Proposition 4.5.1
contains a construction suggested by V.~Kharlamov, which
simplified considerably my original arguments. I should thank also R.
Fintushel for a useful remark about the sign in the product
formula \cite{MST} and B.-L. Wang for commenting his paper
\cite{CW} (which convinced the author that his product formula
implies A4).

\heading\S2. Absolute and relative $\SC$ structures
\endheading
\subheading{2.1. Absolute $\Spinc$ structures} A $\SC$ structure
in a principle $\SO_n$ bundle, $P\to X$, is an isomorphism class
of $\SC_n$-extensions, $S\to P$, of $P$.
The set of $\SC$ structures, $\SC(P)$, has a natural action
$\SC(P)\times H^2(X)\to\SC(P)$, $(\ss,h)\mapsto\ss+h$, which makes
it an affine space over $H^2(X)$.
The projection $\SC_n\to\SO_n\times U_1\to U_1$ associates to $S$
its determinant $U_1$-bundle, $\det S$, with the Chern class
$c_1(S)=c_1(\det S)$, so that $c_1(\ss+h)=c_1(\ss)+2h$.

We simplify the notation writing just $\SC(X)$ instead of
$\SC(P)$, if a principal bundle $P\to X$ is associated with an
obvious vector bundle $E\to X$, for example, with the tangent
bundle of a manifold $X$ (the choice of the euclidian structure in
$E$ is not essential).

\subheading{2.2. The conjugation involution} The conjugate,
$\overline S\to X$, to a principal $\SC_n$ bundle, $S\to X$,
set-theoretically coincides with the latter, but has the conjugate
action of $\SC_n$ (induced by the conjugation automorphism in
$\SC_n$, which covers the direct product automorphism of
$\SO_n\times U_1$, identical on $\SO_n$ and non-identical in
$U_1$).
 The conjugation defines an  involution,
$\conj_P\:\SC(P)\to\SC(P)$, $\ss\mapsto\bar\ss$, such that
$c_1(\bar\ss)=-c_1(\ss)$ and $\overline{\ss+h}=\bar\ss-h$ for any
$\ss\in\SC(P)$ and $h\in H^2(X)$.

\subheading{2.3. Homology interpretation of $\SC$ structures} It
is convenient to identify the set $\SC(P)$ with the coset of the
image of $H^2(X)$ under the monomorphism $\pi_P^*\:H^2(X)\to
H^2(P)$. This makes transparent the nature of the affine structure
in $\SC(P)$. Namely, a $\SC$ extension $F\:S\to P$ can be viewed
as a principal $U_1$-bundle over $P$ since
$\ker(\SCn\to\SO_n)\cong U_1$, and the Chern class $c_1(F)\in
H^2(P)$ defines the correspondence between $\SC$ structures and
those cohomology classes which have a non-trivial restriction
$H^2(P)\to H^2(\SO_n)\cong\Z/2$, $n\ge3$, to a fiber of $P$.

Given $\ss\in\SC(P)\subset H^2(P)$, one can observe that $\pi_P^*(
c_1(\ss))=2\ss$ and that $\overline\ss=-\ss$.

\subheading{2.4. The canonical $\Spinc$ structure of a Lefschetz
fibration} An almost complex structure in a $\SO_{2n}$ bundle
defines the $\SC$-extension associated to the natural homomorphism
$U_n\to\SC_n$. In particular, a symplectic manifold carries a
canonical $\SC$ structure represented by
 {\it the symplectic $\SC$ extension}, $S\to X$.
It is well-known that the total space $X$ of a Lefschetz fibration
$p\:X\to S^2$, carries a compatible  symplectic structure (except
the case of null-homologous fibers of genus $1$, in which $X$ is
still almost complex) which gives the associated canonical $\SC$
structure.

In fact, to define a $\SC$ structure in a vector bundle $E\to X$
it is sufficient to have an almost complex structure  over its
$3$-skeleton, $\Ske_3X$, only (more precisely, $\SC$ structures
can be viewed as equivalence classes of those almost complex
structures over $\Ske_2X$ that can be extended to $\Ske_3X$). This
gives an alternative way to introduce the canonical $\SC$
structure in the Lefschetz fibration, using $U_1\times U_1$
reduction of the tangent bundle $\tau_X$ in the complement of the
critical point set of $p$ determined by the ``vertical'' and the
complementary ``horizontal'' $\SO_2=U_1$ subbundles.

\subheading{2.5. Relative $\Spinc$-structures} In the definition
of a relative structure, in addition to the $\SC_n$-bundle $S\to
X$ considered in 2.1, we fix an isomorphism between the
restriction $S|_A$ and a certain reference principal $\SC$ bundle,
$S_A\to A$. Such a reference bundle appears naturally for example
if $X$ is a $4$-manifold and $A$ is a surface $\SS\subset X$, or a
tubular neighborhood $N$ of $\SS$, or the boundary $\partial N$,
since in these cases $\tau_X|_A$ admits a natural $U_1\times U_1$
reduction and thus, the associated $\SC_4$-extension.

More formally speaking, let $(X,A)$ be a CW-pair, $\pi_P\:P \to X$
a principal $\SO_n$ bundle, $\pi_{P|A}\:P_A\to A$ the restriction
of $\pi_P$ over $A$, and $F_A\:S_A\to P_A$ a $\SC$ extension of
$P_A$. {\it A relative $\Spinc$ extension of $P$} with respect to
$S_A$ is a $\SC$-extension, $S\to X$, $F\:S\to P$, together with
an isomorphism $R\:S|_A\to S_A$, such that $F_A\circ R$ is the
restriction of $F$ to $A$. An isomorphism between relative
$\Spinc$-extensions  $F^{(i)}\:S^{(i)}\to P$, $R^{(i)}\:
S^{(i)}|_A\to S_A$, $i=1,2$, is defined as an isomorphism
$S^{(1)}\to S^{(2)}$ of $\Spinc$ bundles whose restriction over
$A$ commutes with $R^{(1)}$ and $R^{(2)}$. An isomorphism class of
relative $\SC$-extensions is called {\it a relative $\SC$
structure}, and the set of such structures is denoted by
$\Spinc(P,S_A)$, or simply by $\SC(X,A)$ if $P$ and $S_A$ are
evident.

It is straightforward to check that $\Spinc(P ,S_A)$ is an affine
space over $H^2(X,A)$ and the natural forgetful map
$\abs\:\SC(P,S_A) \to \SC(P )$ is affine with respect to the
cohomology forgetful homomorphism $H^2(X,A)\to H^2(X)$.

The conjugation involution defined like in the absolute case
interchanges $\SC(P,S_A)$ with $\SC(P,\overline{S}_A)$. It  is
anti-affine, that is $\overline{\rr+h}=\overline{\rr}-h$ for
$\rr\in\SC(P,S_A)$, $h\in H^2(X,A)$.

\subheading{2.6. Relative $\SC$ structures with respect to
surfaces, $\SS\subset X$, and their evaluation on membranes} Let
$S_\SS\to\SS$ denote the canonical $\SC$ extension defined by the
$U_1\times U_1$ reduction due to the splitting of the tangent
bundle $\tau_X|_\SS$ along $\SS$ into a sum
$\tau_\SS\oplus\nu_\SS$ of the tangent and the normal bundles to
$\SS$.
 Note that the inversion of the orientation of $\SS$
results in the conjugation of the associated canonical $\SC_4$
bundle, $S_{-\SS}=\overline S_\SS$. In particular, the conjugation
involution in this case is $\SC(X,\SS)\to\SC(X,-\SS)$.

Assume that $\rr\in\SC(X,\SS)$. Note that any membrane
$(F,\partial F)\subset(X,\SS)$ defines  trivializations of the
both $\tau_\SS$ and $\nu_\SS$ along $\partial F$, and thus
provides a trivialization of the determinant bundle $\det
S_\SS\cong \tau_\SS \otimes \nu_\SS$. The obstruction class in
$H^2(F,\partial F)$ for extension of this trivialization to the
whole $F$, as it is evaluated on the fundamental class,
$[F,\partial F]$, gives an integer denoted by $\rr[F]$. It is easy
to observe that $\overline\rr[F]=\rr[-F]=-(\rr[F])$.

\subheading{2.7. The canonical relative $\SC$ structures in the
case of symplectic or Lagrangian surface, $\SS\subset X$} Assume
now that $\SS\subset X$ is a symplectic surface with respect to
some symplectic structure $\om$ in $X$ that is $\om|_\SS>0$. Then
we can define the canonical {\it symplectic relative  $\SC$
structure}, $\rr_{X,\SS}\in\SC(X,\SS)$ whose image
$\abs(\rr_{X,\SS})\in\SC(X)$ is the absolute symplectic canonical
$\SC$ structure introduced in section 2.4. Namely, the structure
$\rr_{X,\SS}$ is represented by a $\SC$ extension $S\to P$ of the
principal $\SO_4$ bundle $P\to X$, which arises from an almost
complex structure determined in $\tau_X$ as we fix a Riemannian
metric in $X$ compatible with $\om$. If we choose such a metric
making the surface $\SS$ pseudo-holomorphic (which is always
possible), then the restriction $S|_\SS$ is naturally identified
with the canonical $\SC$-bundle $S_\SS\to\SS$.

Consider now the case of a Lagrangian surface $\SS$ in a
symplectic manifold $X$, in which we can similarly define the
canonical {\it Lagrangian relative $\SC$ structure},
$\rr_{X,\SS}\in\SC(X,\SS)$. One way to do it is to make a
Lagrangian surface symplectic by a perturbation. There is also an
alternative description of this structure (which concerns also the
case of a null-homologous torus in which such a perturbation is
impossible, and makes evident that the choice of a perturbation is
not essential), in which we use the canonical isomorphism between
the two different almost complex structure in $\tau_X|_\SS$: the
one induced from $X$ that is coming from the isomorphism
$\tau_X|_\SS\cong\tau_\SS\otimes\C$, and the one arising from
$U_1\times U_1$ reduction due to the splitting
$\tau_X|_\SS\cong\tau_\SS\oplus\nu_\SS$. It is just
 a special case of the canonical isomorphism
$\xi\otimes\C\cong\xi\oplus\bar\xi$ for a complex bundle $\xi$.
The induced isomorphism of the associated $\SC$ bundles covers an
automorphism of $\tau_X|_\SS$ that can be canonically connected to
the identity by an isotopy. This gives an isomorphism between
these $\SC$ extensions, which defines the Lagrangian relative
$\SC$ structure.

\rk{Remark} Note that the canonical $\SC$ structure
$\rr_{X,\SS}\in\SC(X,\SS)$ of a symplectic pair $(X,\SS)$ is
invariant under the monodromy induced by any symplectic isotopy of
$\SS$ in $X$, whereas any other structure,
$\rr=\rr_{X,\SS}+h\in\SC(X,\SS)$, $h\in H^2(X,\SS)$, is sent by
the monodromy to $\rr_{X,\SS}+f_*(h)$, where $f_*$ is the
cohomology monodromy.

The same concerns Lagrangian surfaces and Lagrangian isotopy.
\endrk

\subheading{2.8. Lefschetz fibrations and their conjugates} A
special case of our interest is  $\SS$ being a fiber of a
Lefschetz fibration $p\:X\to S^2$, (or more generally, a fiber in
a Lefschetz pencil). Such a fiber $\SS\subset X$ is symplectic
with respect to the symplectic form $\om$ supported by the
Lefschetz fibration (or pencil), so there is a canonical structure
$\rr_{X,\SS}\in\SC(X,\SS)$ from section 2.7.

{\it The conjugate Lefschetz fibration} $\overline p\:\overline
X\to \overline{S^2}$ is by definition, set-theoretically the same
as $p$, however, with the opposite orientation chosen in the
base-space $\overline{S^2}=-S^2$ and in the fibers,
$\overline\SS=-\SS$ (so that $\overline X$ itself has the same
orientation as $X$). It is not difficult to observe that
$\overline{\rr_{X,\SS}}= \rr_{\overline X,\overline{\SS}}
\in\SC(X,-\SS)$.

\subheading{2.9. The excision and the homotopy invariance theorems
for $\SC$ structures} The propositions stated below mimic the
standard results for the cohomology and follow automatically from
the latters, since an affine map associated with an isomorphism
must be an affine isomorphism.

\tm{2.9.1. Proposition (excision)} Assume that a CW complex $Z$ is
decomposed into a union of subcomplexes, $Z=X\cup Y$, $A=X\cap Y$.
 Consider a principal $\SO_n$ bundle $P_Z\to Z$ and let $P_X$, $P_Y$,
$P_A$ denote its restrictions over $X$, $Y$, and $A$ respectively.
Fix a $\SC$ extension, $F_Y\:S_Y\to P_Y$ and let $F_A\:S_A\to P_A$
denote its restriction over $A$. Then the restriction map
$$ \SC(P_Z,S_Y)\to\SC(P_X,S_A)$$
is an isomorphism of affine spaces agreeing with the isomorphism $H^2(Z,Y)\cong
H^2(X,A)$. \qed\endtm

\tm{2.9.2. Proposition (homotopy invariance)} Assume that  $\SS$
is a deformation retract of $N\subset X$. Let $F_N\:S_N\to P_N$ be
a $\SC$ extension and $S_\SS=S_N|_\SS$. Then the restriction map
$\SC(P_X,S_N)\to\SC(P_X,S_\SS)$ is an isomorphism of the affine
spaces agreeing with the isomorphism $H^2(X,N)\cong H^2(X,\SS)$.
\qed\endtm

\tm{2.9.3. Corollary} Let $\SS\subset X$ be a surface in a
$4$-manifold, $N\subset X$ its compact tubular neighborhood,
$M=\partial N$, and $X^\circ=X\smallsetminus\Int(N)$. Then we have
canonical affine isomorphisms
$\SC(X,\SS)\cong\SC(X,N)\cong\SC(X^\circ,M)$. \qed\endtm

\subheading{2.10. Connected sums and blowing up of  $\SC$
structures} Definitions of the connected sum and the blowup
operations are obvious and well-known for the absolute $\SC$
structures. They can be easily extended to relative $\SC$
structures, as well. For given $n$-manifolds $X_\pm$, with
codimension $2$ submanifolds $\SS_\pm$, and
$\rr_\pm\in\SC(X_\pm,\SS_\pm)$, we obtain the connected sum
$\rr_+\#\rr_-\in\SC(X_+\# X_-,\SS_+\#\SS_-)$, where
$\SS_+\#\SS_-\subset X_+\# X_-$ is the internal connected sum of
$(X_\pm,\SS_\pm)$.

Let $\SS\subset X$ be a surface in a four-manifold, and
$\hat\SS\subset\hat X$ its proper image after blowing up $X$ at a
point of $\SS$, that is $\SS\#\Cp1\subset X\#(-\Cp2)$. For
$\rr\in\SC(X,\SS)$, we define $\hat\rr\in\SC(\hat X,\hat\SS)$ as
$\hat\rr=\rr\#\rr_{-1}$, where $\rr_{-1}\in\SC(-\Cp2,\Cp1)$ is the
unique structure such that $c_1(\abs(\rr_{-1}))=-1$.

\subheading{2.11. Fiber sums of relative $\SC$ structures}  Let
$\SS$ be a closed oriented surface of genus $g>1$. We say that $X$
is {\it a $\SS$-marked $4$-manifold}, if there is a fixed smooth
embedding $f\:\SS\to X$ endowed with a normal framing of $f(\SS)$
(in particular, $\SS^2=0$). To simplify the notation, we will be
writing $\SC(X,\SS)$ rather then $\SC(X,f(\SS))$.

Given $\SS$-marked $4$-manifolds $X_\pm$, consider their fiber sum
$X=X_+\#_\SS X_-=X_+^\circ\cup_f X_-^\circ$, where
$X_\pm^\circ=X_\pm\smallsetminus\Int(N_\pm)$ (the complements of
the tubular neighborhoods of $\SS$), and the gluing diffeomorphism
$f\:\partial N_+\to\partial N_-$ is naturally determined by the
trivialization of $N_\pm\to\SS$ respecting the framings, so that
$\partial N_+$ and $-\partial N_-$ are identified with
$M=\SS\times S^1$.
 Note that $X$ has an induced structure of $\SS$-marked $4$-manifold, since
$\SS_t=\SS\times t\subset\SS\times S^1$ has a natural normal
framing.

Operations $\rr_+\#_\SS\rr_-\in\SC(X)$ and
$\rr_+\vee\rr_-\in\SC(X,\SS)$ for $\rr_\pm\in\SC(X_\pm,\SS)$ are
the compositions of the isomorphism
$$
\SC(X_+,\SS)\times \SC( X_-,\SS)\cong \SC(X^\circ_+,M)\times
\SC(X^\circ_-,M)\cong\SC(X,M)
$$
with the forgetful maps $\SC(X,M)\to\SC(X)$ and
$\SC(X,M)\to\SC(X,\SS)$.

Given $\ss_\pm\in\SC(X_\pm,\ss_\SS)$ we denote by
$\ss_-\#_\SS\ss_+$ a subset of $\SC(X)$ consisting of the
structures $\rr_-\#_\SS\rr_+$ for all
$\rr_\pm\in\abs^{-1}_{X_\pm,\SS}(\ss_\pm)$. It is not difficult to
check that the set $\ss_-\#_\SS\ss_+$ is affine with respect to
the subgroup $\Delta_M\subset H^2(X)$, which is the image of
$H^1(\SS)$ under the product of the homomorphism $q^*\:H^1(\SS)\to
H^1(M)$ induced by the projection $q\:M\cong\SS\times S^1\to\SS$
and the boundary map $\delta_M\: H^1(M)\to H^2(X)$.

One can also interpret $\ss_-\#_\SS\ss_+$ as a set consisting of those $\ss\in\SC(X)$
which have $d(\ss)=d(\ss_+)+d(\ss_-)$ and whose restriction to $X_\pm^\circ$ coincides
with that of $\ss_\pm$.

\subheading{2.12. The natural properties of the operations with
the relative $\SC$ structures} It is not difficult to check that
the operations introduced in 2.11 satisfy the following natural
properties
$$\align
\rr_1\vee\rr_2=\rr_2\vee\rr_1,\ \ &\text{and thus }
\ \rr_1\#_\SS\rr_2=\rr_2\#_\SS\rr_1 \\
(\rr_1\vee\rr_2)\vee\rr_3=\rr_1\vee(\rr_2\vee\rr_3),\ \ &\text{and
thus }\ (\rr_1\vee\rr_2)\#_\SS\rr_3=\rr_1\#_\SS(\rr_2\vee\rr_3) \\
\overline{\rr_1\vee\rr_2}=\overline{\rr}_1\vee\overline{\rr}_2 ,\
\ &\text{and thus }\ \overline{\rr_1\#_\SS\rr_2}=
\overline{\rr}_1\#_{-\SS}\overline{\rr}_2\\
\endalign$$
where the equalities mean that the obvious diffeomorphisms
$$\align(X_1\#_\SS X_2,\SS)&\cong(X_2\#_\SS X_1,\SS)\\
((X_1\#_\SS X_2)\#_\SS X_3,\SS)&\cong (X_1\#_\SS(X_2\#_\SS
X_3),\SS)\\X_1\#_\SS X_2&\cong X_1\#_{-\SS}X_2\endalign$$ send one
of the corresponding $\SC$ structures to the other.

Some ambiguity in the notion of ``the obvious diffeomorphism'',
related in particular to the ambiguity in $\SS$-marking of the
fiber sums, turns out to be not essential.
It is also straightforward to check the following

\tm{Proposition 2.12.1} Assume that $(X_i,\SS)$ are symplectic
pairs, $i=1,2$, $X=X_1\#_\SS X_2$, and $\rr_i\in\SC(X_i,\SS)$ are
the canonical relative $\SC$ structures. Then
$\rr_1\vee\rr_2\in\SC(X,\SS)$ is also the canonical $\SC$
structure of the symplectic pair $(X,\SS)$. In particular,
$\rr_1\#_\SS\rr_2\in\SC(X)$ is the canonical symplectic $\SC$
structure of $X$. \qed\endtm

\heading \S3. The basic properties of the absolute SW and OS
invariants
\endheading

\subheading{3.1. The axioms} Axioms A1, A3 and A4 below are
essential for the definition of $S_{X,\SS}$, for showing its
independence from $W$, whereas axioms A2 and A5 are required only
for proving the corresponding properties of $S_{X,\SS}$, namely,
(3) and (7) in Theorem 1.4.1.
 Unless it is stated otherwise,
we suppose in this section that all the closed $4$-manifolds below
have $b_2^+>1$ (in the case $b_2^+=1$ the formulations are
similar, but require a bit more care).

{\bf A1.} {\it Finiteness.} The set of the basic $\SC$ structures
$B_X=\{\ss\in\SC(X)\,|\, S_X(\ss)\ne0\}$ is finite for any $X$.

{\bf A2.} {\it Conjugation symmetry}.
$S_X\circ\conj_X=\pm S_X$, where $\conj_X$ is the conjugation
involution in $\SC(X)$.

 {\bf A3}. {\it Lefschetz
normalization.} Assume that $X\to S^2$ is a relatively minimal
Lefschetz fibration whose fiber $\SS\subset X$ has genus $g>1$.
Let $\ss_X\in\SC(X)$ denote the canonical $\SC$ structure.
 Then \roster\item
$S_X(\ss_X)=1$, if $X$ is endowed with the canonical homology orientation (with
respect to a symplectic structure supporting the Lefschetz fibration).
\item
$\ss_X$ is the only basic structure $\ss\in\SC(X)$ satisfying the
condition $c_1(\ss)[\SS]=\chi(\SS)+\SS^2$; \item for any fiber sum
$X\#_\SS Y=X^\circ\cup Y^\circ$ with a $\SS$-marked $4$-manifold,
the restriction $\ss|_{X^\circ}\in\SC(X^\circ)$ of any basic $\SC$
structure $\ss\in\SC(X\#_\SS Y,\ss_\SS)$,  coincides with the
restriction $\ss_X|_{X^\circ}\in\SC(X^\circ)$.
\endroster

{\bf A4.} {\it Product formula.}  Let $X={X}_-\#_\SS
X_+=X_-^\circ\cup X_+^\circ$ be a fiber sum like in 2.11, with a
fiber $\SS$ of genus $g>1$. Choose $\ss_\pm\in\SC(X_\pm,\ss_\SS)$
and let $\s_\pm\in H^2(X_\pm)$ denote the Poincare dual class to
$\SS\subset X_\pm$.
 Then
$$\sum_{k\in\Z}S_{X_-}(\ss_--k\s_-)S_{X_+}(\ss_++k\s_+)=\sum_{\ss\in\ss_-\#_\SS\ss_+}
S_X(\ss).$$

{\bf A5.} {\it Adjunction inequality:}
 $-\chi(\SS)\ge\SS^2+|c_1(\ss)[\SS]|$,
for any $\ss\in B_X$ and an essential surface $\SS\subset X$ of
genus $g>0$, with $\SS^2\ge0$.

Combining properties A3(1) and A4 and taking into account the
remark about $\ss_-\#_\SS\ss_+$ and $\Delta_M$ in the end of
section 2.11, we obtain the following

\tm{3.1.1. Corollary} If the summands, $X_-$, involved into a
fiber sum in A$4$ is a Lefschetz fibration with a fiber $\SS$, then
for any $\ss_+\in\SC(X_+,\ss_\SS)$
$$
S_{X_+}(\ss_+)=\sum_{\ss\in\ss_-\#_\SS\ss_+}
S_X(\ss)=\sum_{\rr_+\in\abs^{-1}(\ss_+)} S_X(\rr_+\#_\SS\rr_-)
=\sum_{h\in\Delta_M} S_X(\ss+h)
$$
where $\ss_-\in\SC(X_-)$, $\rr_-\in\SC(X_-,\SS)$ are the canonical
absolute and relative $\SC$ structure of the Lefschetz fibration.
In the last sum, $\ss$ is any fixed element of
$\ss_-\#_\SS\ss_+\subset\SC(X)$.
\endtm

\subheading{3.2. Properties A1, A2, and A5}

{\bf A1.} The finiteness is a fundamental well-known property of
SW invariants, which holds as well for OS invariants, see
\cite{OS3}, Theorem 3.3.

{\bf A2.} It is also a well-known property. In fact, a
set-theoretic identification of the conjugate $\SC$-bundles $S$
and $\bar S$ gives a point-wise correspondence (possibly
alternating the orientations) between the solutions spaces to the
SW equations associated with $S$ and $\bar S$. For the case of OS
invariants, see \cite{OS3}, Theorem 3.5.

{\bf A5.} We formulated the simplest classical version of the
adjunction inequality. It can be found in a more general form
(including the case of $\SS^2<0$) in \cite{OS2}, Theorems 1.1--1.7
for SW invariants, and in \cite{OS3}, Theorem 1.4. for OS
invariants.

\subheading{3.3. Lefschetz normalization properties A3(1)--(3)}

{\bf A3(1).}  For SW invariants this property is proved by Taubes
\cite{T1} (the Main Theorem). The case of OS invariants was
considered in \cite{OS4}, Theorem 5.1. For this part of A3 the
minimality condition is not required.

{\bf A3(2).} By Taubes' result \cite{T2}, Theorem 2,
$|c_1(\ss)\circ[\om]|\le|c_1(\ss_X)\circ[\om]|$ for any SW basic
structure $\ss\in\SC(X)$, with the equality only for $\ss=\ss_X$
and $\ss=\bar\ss_X$ (here $\circ$ denotes the pairing in
$H^2(X)$). According to \cite{OS4, Theorem 1.1}, the same holds
for OS basic structures.
Gompf's construction produces a symplectic form in a Lefschetz
fibration $p:X\to S^2$  as a small perturbation of
$\om=p^*(\om_{S^2})+t\eta$, where $\om_{S^2}$ is the area form in
$S^2$, $\eta$ is a closed $2$-form in $X$ having positive
restrictions to the fibers of $p$ at every point, and $0<t<\!<1$.
Observing that $[\SS]$ is dual to a properly normalized $2$-form
$p^*\om_{S^2}$ and, thus,
 $c_1(\ss)\circ p^*\om_{S^2}=c_1(\ss)[\SS]$, we can deduce
letting $t\to0$ that $|c_1(\ss)[\SS]|\le|c_1(\ss_X)[\SS]|$. Using
that $[\eta]\in H^2(X)$ can be any class with $[\eta][\SS]>0$, we
can also conclude that the equality
$|c_1(\ss)[\SS]|=|c_1(\ss_X)[\SS]|$ may hold only in the case of
$\ss=\ss_X+n\sigma$, or $\ss=\bar\ss_X+n\sigma$, where $\sigma\in
H^2(X)$ is dual to $[\SS]$. But the symplectic manifolds are of
the simple type \cite{T4}, Theorem 02(6), which implies that
$\ss_X+n\sigma$ (and similarly $\bar\ss_X+n\sigma$) cannot be
basic for $n\ne0$, since $[c_1^2(\ss_X+n\sigma)
-c_1^2(\ss_X)]=2n\chi(\SS)\ne0$ in case of $g(\SS)>1$.

{\bf A3(3) for OS invariants.}
 It is proved by the arguments in Lemma 5.7 from
 \cite{OS4} for OS invariants. The same scheme of the proof works for
SW invariants as well, so we will briefly review it (sending a
reader to \cite{OS4} for the notation and details).

The first step is to observe that the canonical structure
$\ss_X\in\SC(X)$ is the only one satisfying the adjunction
inequality with respect to a certain family of surfaces $F\subset
X$. For these surfaces $F^2<0$, and so in principle the inequality
may fail for a basic structure $\ss\in\SC(X)$, but in this case
there is another basic $\SC$ structure ${\ss'}=\ss+f$, where $f\in
H^2(X)$ is Poincare dual to $[F]$, and there exists $\xi\in\A_F$
such that $\Phi_{X,{\ss'}}(\xi x)=\Phi_{X,\ss}(x)$ for any
$x\in\A_X$ (the action of $\xi$ on $x$ means the action of the
image of $\xi$ under the inclusion map $\A_F\to\A_X$).
One can notice next that the construction of surfaces $F$ in
\cite{OS4} yields a natural epimorphism $H_1(\SS)\to H_1(F)$
commuting with the inclusion homomorphisms from $H_1(\SS)$ and
$H_1(F)$ to $H_1(X)$, and so we may assume that $\xi\in\A_\SS$.

The second key observation is triviality of the action of $\A_\SS$
in $HF^+(M,\ttt)$, where $M=\SS\times S^1$ and $\ttt=\ss|_M$ is
the canonical structure induced from $\ss_\SS\in\SC(\SS)$ by the
projection $M\to\SS$. This triviality is deduced in \cite{OS4} as
an immediate corollary of the isomorphism $HF^+(M,\ttt)\cong\Z$.

The third ingredient of the proof is the relation between the
invariants $F^\mix_{X^\circledcirc,{\ss'}^\circledcirc}$ and
$\Phi_{X,\ss'}$ (the latter is dual to $\OS_X(\ss')$ in our notation),
where $X^\circledcirc=X\smallsetminus(\Int N\cup \Int B^4)$ is the
complement of a tubular neighborhood $N$ of $\SS$ and a ball
$B^4\subset X$ disjoint from $N$. $X^\circledcirc$ is viewed as a
cobordism from $S^3$ to $M=\partial N$, so that
$F^\mix_{X^\circledcirc,{\ss'}^\circledcirc}$ takes values in
$HF^+(M,\ttt)\cong\Z$. More precisely, $\Phi_{X,{\ss'}}(\xi x)$
and thus $\Phi_{X,{\ss}}(x)$, or equivalently, $\OS_X(\ss)$,
vanishes as it is the homogeneous part of
$F^\mix_{X^\circledcirc,{\ss'}^\circledcirc}(\Theta_-\otimes(\xi
x))=\xi
F^\mix_{X^\circledcirc,{\ss'}^\circledcirc}(\Theta_-\otimes x)$,
where $\Theta_-\in HF^-(S^3)$ is the generator in the upper
dimension, and ${\ss'}^\circledcirc={\ss'}|_{X^\circledcirc}$ (see
the proof of Lemma 5.6 of \cite{OS4}).
This contradicts to the assumption that $\ss$ (and thus $\ss'$) is
a basic structure.

Applying these arguments to $X\#_\SS Y$, we conclude similarly
that if $\OS_{X\#_\SS Y}(\ss)\ne0$, for $\ss\in\SC(X\#_\SS Y)$
such that $\ss|_\SS=\ss_\SS$, then $\ss|_{X^\circ}$ is canonical,
since otherwise $F^\mix_{X^\circ\#_\SS Y,{\ss'}^\circ}$ and, thus,
$\OS_{X\#_\SS Y}(\ss)$ vanish.\qed

{\bf A3(3) for SW invariants.}
The first step for SW invariants is like for OS invariants, since
the generalized adjunction inequalities look similar in the both
theories (cf. \cite{OS1}, \cite{OS2} and \cite{OS4}). Next, the
Seiberg-Witten-Floer homology groups $HF^{\SW}_*(M,\ttt)\cong\Z$,
(see \cite{MW}, Theorem 1.7), and so the action of $\A_\SS$
considered in \cite{CW} is trivial on this group for the same
reason as in the case of the OS invariants.

The final step goes also like in the OS-theory, but instead of
$F^\mix_{X^\circledcirc,{\ss'}^\circledcirc}(\Theta_-\otimes x)$ we consider the function
$\phi^{\SW}_{X^\circ}({\ss'}^\circ,x^\circ)$ from \cite{CW}, where
$X^\circ=X\smallsetminus \Int N$, $x^\circ=x|_{X^\circ}$ and
${\ss'_X}^\circ=\ss'|_{X^\circ}$.
To deduce vanishing of $\SW_{X,\ss'}(\xi x)$ we can use the product formula,
in Theorem 1.2 of \cite{CW}, which implies for $X\#_\SS
Y=X^\circ\cup Y^\circ$ that
$$\SW_{X\#_\SS Y,\ss'}(\xi x\otimes y)=
\la[u]\pi_1(\phi^{\SW}_{X^\circ}({\ss'_X}^\circ,\xi
x)),\pi_2(\phi^{\SW}_{Y^\circ}({\ss'_Y}^\circ,y))\ra
$$
and we can conclude that the product vanishes, because the first
factor vanishes. \qed

\subheading{3.4. The product formula A4}
 The well known product formula
\cite{MST}, Theorem 3.1, concerns the version of SW invariants
corresponding to the case of $R_X=\Z[U]$, which is not as general
as $R_X=\A_X^*$, although the author supposes that the same
arguments without essential changes can be as well used in the most general
case. In any case, after \cite{MST}, much more general gluing
formulae were established, see for instance Theorem 1.2 in
\cite{CW}, which concerns $4$-manifolds with an arbitrary boundary
and contains A4 as a corollary.

In the case of OS invariants, A4 can be derived from the product
formula \cite{OS3}, Theorem 3.4, applied to the fiber sums,
although it may look not so obvious as in the case of SW
invariants. To clarify it, we give some comments, which are
basically extracted from \cite{OS3} and \cite{OS4}.

Puncturing a fiber sum, $X=X_-\#_\SS X_+=X_-^\circ\cup X_+^\circ$,
at a pair of points, we obtain $X^\circledcirc=X_-^{\circledcirc}
\cup X_+^\circledcirc$, where
$X_\pm^{\circledcirc}=X_\pm^\circ\smallsetminus B^4$.
$X^\circledcirc$ can be viewed as a product of cobordisms
$X_-^{\circledcirc}\:S^3\to M\cong S^1\times \SS$ and
$X_+^\circledcirc\: M\to S^3$. The product formula \cite{OS3},
Theorem 3.4, says
$$
[F^+_{X_+^{\circledcirc},\ss_+^{\circledcirc}}(
F^\mix_{X_-^{\circledcirc} , \ss_-^{\circledcirc}}(\theta_-
\otimes x_-)\otimes x_+)]_0 = \sum_{\ss|_{X_\pm^{\circledcirc}}=
\ss_\pm^{\circledcirc}} \Phi_{X,\ss}(x_-\otimes x_+)
$$
where $\ss_\pm^{\circledcirc}\in\SC({X_\pm^{\circledcirc}})$,
$\ss_\pm^{\circledcirc}|_M=\ttt$, and $[x]_0\in HF^+_0(S^3)$
denotes the $0$-dimensional component of $x\in HF^+(S^3)$.

We assume here that $\SC$ structure $\ss_\pm^{\circledcirc}$ is
induced in $X^\circledcirc_\pm$ from
$\ss_\pm\in\SC(X_\pm,\ss_\SS)$, and, thus, its restriction, $\ttt$
is the canonical $\SC$ structure determined by the $\SO_2$
reduction of $\tau_M$. This implies, in particular, that
$HF^+(M,\ttt)\cong\Z$ \cite{OS4}, Lemma 5.5.

The duality between
$F^+_{X_+^{\circledcirc},\ss_+^{\circledcirc}}\:HF^+(M,\ttt)\to
HF^+(S^3)$ and
$F^-_{X_+^{\circledcirc},\ss_+^{\circledcirc}}\:HF^-(S^3)\cong
HF^-(-S^3)\to HF^-(-M,\ttt)\cong HF^-(M,\ttt)$ (see Theorem 3.5 in
\cite{OS3}) implies that
$$[F^+_{X_+^{\circledcirc},\ss_+^{\circledcirc}} (1\otimes x_+)]_0=
F^\mix_{X_+^{\circledcirc} , \ss_+^{\circledcirc}}(\theta_-\otimes
x_+),$$ which gives
$$
F^\mix_{X_+^{\circledcirc},\ss_+^{\circledcirc}}(\theta_-\otimes
x_+) F^\mix_{X_-^{\circledcirc}, \ss_-^{\circledcirc}}(\theta_-
\otimes x_-) = \sum_{\ss|_{X_\pm^{\circledcirc}}=
\ss_\pm^{\circledcirc}} \Phi_{X,\ss}(x_-\otimes x_+)
$$

On the other hand, applying the product formula \cite{OS3}, to
$X_\pm^\circledcirc=X_\pm^\circledcirc\cup
N_\pm^\circledcirc(\SS)$, viewed as a product cobordism of
$X_\pm^\circledcirc\:S^3\to M$ and $N_\pm^\circledcirc\: M\to
S^3$, we obtain
$$
F^\mix_{X^\circledcirc_\pm,\ss_\pm^\circledcirc}(\theta_-\otimes
x_\pm)=
\sum_{\ss|_{X_\pm^\circledcirc}=\ss^\circledcirc}\Phi_{X_\pm,\ss_\pm}(x_\pm)
$$
using that the second cobordism induces an isomorphism from
$HF^+(M,\ttt)\cong\Z$ to $HF^+_0(S^3)$ (see \cite{OS4}, Theorem
5.3).

The structures $\ss\in\SC(X_\pm)$ in the latter sum differ just by
multiples of the class $\s\in H^2(X_\pm)$ Poincare-dual to
$[\SS]$, and, thus, have distinct degrees, $d(\ss)$, since
$d(\ss+n\s)=d(\ss)+n\chi(\SS)$. In the other words, the latter
formula is a decomposition of
$F^\mix_{X_\pm^\circledcirc,\ss^\circledcirc_\pm}$ into a sum of
its homogeneous components (this idea was used in the proof of
Lemma 5.6 in \cite{OS4}).
Passing from $\Phi_{X,\ss}$ to the dual $\OS_X(\ss)$ and comparing
the components of the same degree, we obtain A4.

\comment The complement $X^\circ=X\smallsetminus N(\SS)$ of a
tubular neighborhood of $\SS\subset X$ punctured at a point gives
a cobordism from $S^3$ to $M=\partial X^\circ\cong\SS\times S^1$,
which induces homomorphism \cite{OS3}
$$F^\mix_{X^\circ\smallsetminus B^4,\ss}\:HF^-(S^3) \otimes \A_{X^\circ}\to
HF^+(M,\ttt)$$

\endcomment

\heading \S4. Proof of Theorems 1.4.1 \endheading

\subheading{4.1. Independence of the choice of $W$} Consider a
pair of Lefschetz fibrations, $W_i\to S^2$, $i=1,2$, with a fiber
$\SS$ such that $H_1(W_i)=0$, and denote by $\rr_i\in\SC(W_i,\SS)$
the canonical structures. Let $Y_i=X\#_\SS W_i$ and $W=W_1\#_\SS
W_2$, then $Z=X\#_\SS W\cong X\#_\SS W_1\#_\SS W_2\cong Y_1\#_\SS
W_2$.
\tm{4.1.1. Proposition} For any $\rr\in\SC(X,\SS)$ we have
$$S_{Y_1}(\rr\#_\SS\rr_1)=S_{Z}(\rr\#_\SS\rr_1\#_\SS\rr_2)=S_{Y_2}(\rr\#_\SS\rr_2)$$
\endtm

\demo{Proof} Since the two equalities are analogous, we prove only
the first one. Let $W_i^\circ=W_i\smallsetminus N_i$, $i=1,2$,
denote the complements of an open tubular neighborhood $N_i$ of a
fiber $\SS\subset W_i$, and $Y_i^\circ=X\#_\SS W_i^\circ$,
$W^\circ=W_1^\circ\#_\SS W_2$.
From Corollary 3.1.1 applied to a fiber sum of $Y_1\#_\SS W_2$ we
obtain
$$S_{Y_1}(\rr\#_\SS\rr_1)=\sum_{h\in\Delta_M}S_{Z} ((\rr\#_\SS\rr_1\#_\SS
\rr_2) + h).
$$
where $M=\partial W_2^\circ$ and $\Delta_M\subset H^2(Z)$ is the
image of $H^1(\SS)$ under the composition $\delta_M\circ q^*\:
H^1(\SS)\to H^1(M) \to H^2(Z)$, like in 2.11.

Note that the sum in the above formula has only one non-vanishing
term corresponding to $h=0$, because the restriction of any basic
structure $(\rr\#_\SS\rr_1\#_\SS \rr_2)+h$ to $W$ should coincide
with $\rr_1\#_\SS\rr_2|_{W^\circ}$ according to A3(3) and 2.12.1.
On the other hand, for $h\ne0$ it does not coincide, because of
the following observation.

\tm{4.1.2. Lemma} The following composition is injective
$$\CD H^1(\SS)@>q^*>> H^1(M)@>\delta_M>> H^2(Z)@>>> H^2(W^\circ)\endCD$$
(the last map here is the inclusion homomorphism).
\endtm

\demo{Proof} The Poincare dual to these homomorphisms are the
homomorphisms
$$H_1(\SS)\to H_2(M)\to
H_2(Z)\to H_2(W^\circ,\partial W^\circ)$$ sending
 $h_1\in H_1(\SS)$ to the image of $h_2=h_1\times[S^1]\in H_2(\SS\times S^1)\cong
H_2(M)$ under the inclusion map $H_2(M)\to H_2(W^\circ)$ composed
with the relativization map $H_2(W^\circ)\to H_2(W^\circ,\partial
W^\circ)$. The condition that $H_1(W_i)=H_1(W_i^\circ)=0$ allows
to find a cycle in $H_2(W^\circ)$ having non-vanishing
intersection index with $h_2$ in $W$, if $h_1\ne0$, thus proving
non-vanishing of the image of $[h_2]$ in $H_2(W^\circ,\partial
W^\circ)$. \qed\enddemo
\enddemo

\subheading{4.2. Proof of Properties (1)--(5) in Theorems 1.4.1}

{\bf (1)} This property is just  A1 applied to $X\#_\SS W$.

{\bf (2)} This holds by definition
 of the invariants in the case of $\SS^2>0$.

{\bf (3)} Note that the connected sum $X\#_\SS W$ is the same as
the sum $X\#_{-\SS}\overline W$, where $\overline W$ is the
conjugate to $W$ Lefschetz fibration. Axiom $A2$ implies that
$S_{X,\SS}(\rr)=S_{X\#_\SS W}(\rr\#_\SS\rr_{W,\SS})$ is equal to
$\pm S_{X\#_\SS W }(\overline{\rr\#_\SS\rr_{W,\SS}})$ where the
conjugate $\SC$ structure $\overline{\rr\#_\SS\rr_{W,\SS}}$ equals
to $\overline{\rr}\#_{-\SS}\overline{\rr_{W,\SS}}$ as follows from
2.12, and $\overline{\rr_{W,\SS}}=\rr_{\overline W,\overline{
\SS}}$, as
 remarked in 2.8.
On the other hand, using $\overline W$ to evaluate
$S_{X,-\SS}(\overline\rr)$, we obtain
$S_{X,-\SS}(\overline\rr)=S_{X\#_{-\SS}\overline W} (\overline\rr
\#_{-\SS} \rr_{\overline W,\overline\SS})$ that is $\pm
S_{X,\SS}(\rr)$.

{\bf (4)} It follows immediately from Proposition 2.12.1 and
A3(1).

{\bf (5)} It follows from Corollary 3.1.1 applied to the fiber sum
$X\#_\SS W$, namely
$$S_X(\ss)=\sum_{\rr\in\abs_{X,\SS}^{-1}(\ss)}S_{X\#_\SS W}(\rr\#_\SS\rr_{W,\SS})
=\sum_{\rr\in\abs_{X,\SS}^{-1}(\ss)}S_{X,\SS}(\rr),
$$
where $\rr_{W,\SS}\in\SC(W,\SS)$ is the canonical relative $\SC$
structure of a Lefschetz fibration and
$\abs_{X,\SS}\:\SC(X,\SS)\to\SC(X)$ the forgetful map. \qed

\subheading{4.3. Proof of the product formula (6)} Consider a
fiber sum $X=X_+\#_\SS X_-$ and Lefschetz fibrations $W_\pm\to
S^2$ with a fiber $\SS$ and $H^1(W_\pm)=0$. Put $Y_\pm=X_\pm\#_\SS
W_\pm$, $W=W_+\#_\SS W_-$, and $Z=X\#_\SS W\cong Y_+\#_\SS Y_-$.

Choose a pair of relative structure $\rr_\pm\in\SC(X_\pm,\SS)$, denote by
$\rr_{W_\pm,\SS}\in\SC(W_\pm,\SS)$ the canonical $\SC$ structures of Lefschetz
fibrations in $W_\pm$ and let $\ss_\pm=\rr_\pm\#_\SS\rr_{W_\pm,\SS} \in \SC(Y_\pm)$,
$\ss=\rr_+\#_\SS\rr_{W_+,\SS}\#_\SS\rr_-\#_\SS\rr_{W_-,\SS}\in\SC(Z)$.
 By Proposition 2.12.1,
$\rr_{W,\SS}=\rr_{W_-,\SS}\vee\rr_{W_+,\SS}\in\SC(W,\SS)$
 is the canonical relative $\SC$ structure of the Lefschetz fibration in $W$
and $\ss=(\rr_+\vee\rr_-)\#_\SS\rr_{W,\SS}$.

The product formula A4 applied to the fiber sum $Z=Y_+\#_\SS Y_-$
reads
$$\sum_{k\in\Z}S_{Y_+}(\ss_++k\s_+)S_{Y_-}(\ss_--k\s_-)=
\sum_{\ss'\in\ss_+\#\ss_-}
 S_{Z}(\ss')
$$
 The  sum in the right-hand side contains only one
term $S_Z(\ss)$ which follows from the arguments analogous to
those in 4.1. Finally, we observe that $S_{X_\pm,\SS}(\rr_\pm\pm
k\s_\pm)=S_{Y_\pm}(\ss_\pm\pm k\s_\pm)$,
$S_{X,\SS}(\rr_+\vee\rr_-)=S_{Z}(\ss)$.

Equivalence of the alternative formulations of the product formula
in 1.4.1(6) follows from that $\rr_+'\vee\rr_-'=\rr_+\vee\rr_-$ if
and only if $\rr_\pm'=\rr_\pm\pm k\s_\pm$ for some $k\in\Z$.
 \qed

\subheading{4.4. Proof of the Adjunction inequality (7)} The idea
of the proof is to find an appropriate Lefschetz fibration $W\to
S^2$ with a fiber $\SS_W\cong\SS$ having a membrane, $F_W\subset
W$ whose boundary, $\partial F_W$ matches with the boundary
$\partial F$ of membrane $F\subset X$. Then after gluing $F$ and
$F_W$ we can get a closed surface $\hat F\subset X\#_\SS W$, which
will be oriented if the orientations of $\partial F$ and $\partial
F_W$ match.

More precisely, we should glue the complements $X^\circ$ and $W^\circ$ of tubular
neighborhoods $N\subset X$ of $\SS$ and $N_W\subset W$ of $\SS_W$ so that $F\cap
\partial N$ is glued to $F_W\cap\partial N_W$. It is not difficult to see that
connectedness of $\SS\smallsetminus \partial F$ guarantees that we can find such a
gluing map $\partial N\to\partial N_W$.

Finally, we want to make use of the adjunction inequality A5 for
$\hat F$. This requires $\hat F^2=F^2+ F_W^2\ge0$, which holds if
we can find $F_W$ with a sufficiently big self-intersection index.
If we choose $F_W$ so that $\rr_{W,\SS}[F_W]=0$ for the canonical
structure $\rr_{W,\SS}$ of the Lefschetz fibration, then
$c_1(\rr\#\rr_{W,\SS})[\hat F]=\rr[F]+\rr_{W,\SS}[F_W]=\rr[F]$ and
thus
$$-\chi(\hat F)=-\chi(F)-\chi(F_W)\ge F^2+F_W^2+|\rr[F]|$$
which gives (7) of Theorem 1.4.1 provided $F_W^2=-\chi(F_W)$. So,
we reduced the problem to constructing the following example of
$W$ and $F_W$.

\tm{4.4.1. Proposition} Let $\SS$ be a surface of genus $g\ge1$,
 $L\subset\SS$ an oriented curve (possibly multi-component) with the
connected complement $\SS\smallsetminus L$, and $n\in\N$. Then there exists a
relatively minimal Lefschetz fibration $p\: W\to\Cp1$, with $H_1(W)=0$, whose fiber,
$\SS$, has a membrane, $F_W\subset W$, such that \roster\item $\partial F_W=L$ (as an
oriented curve),
\item
$F_W^2=-\chi(F_W)$,
\item $\rr_{W,\SS}[F_W]=0$,
where $\rr_{W,\SS}\in\SC(W,\SS)$ is the canonical $\SC$ structure,
\item $-\chi(F_W)>n$.
\endroster
\endtm

\subheading{4.5. Real Lefschetz fibrations} We will construct a
complex algebraic Lefschetz fibration $p\:W\to\Cp2$ endowed with a
{\it real structure}, that is just an anti-holomorphic involution
({\it the complex conjugation in $W$}), $c\:W\to W$, which
commutes with $p$ and the complex conjugation in $\Cp1$.
The {\it real locus}, $\R W$, of $W$ is the fixed point set of
$c$. For a real fiber, $\SS=p^{-1}(b)$, $b\in\Rp1$, we let
$\R\SS=\SS\cap\R W$. A membrane $F_W$ in our example will be the
closure of a properly chosen connected component of $\R
W\smallsetminus\R\SS$ bounded by $\R\SS$. Such a choice guarantees
the condition (2) of the Proposition 4.4.1, since the tangent
bundle to
 $\R W$ is anti-isomorphic to the normal bundle via the operator
$J\:\tau_X\to\tau_X$ of the complex structure.
 The condition (3) follows from that the real determinant gives a
 section trivializing the complex determinant bundle
(and thus the associated $\SC$ determinant), provided
$F_W\subset\R W$ is orientable.

Note furthermore that the pairs $(\SS,L)$ are classified up to
homeomorphism respecting the orientations of $\SS$ and $L$ just by
the genus $g$ and the number of components, $r\le g$, of $L$,
under our assumption that $\SS\smallsetminus L$ is connected. So,
we will achieve $(\SS,\R\SS)\cong (\SS,L)$, that is the condition
(1), if $\R\SS$ does not divide $\SS$ into halves and has the
required number $r$ of the components.
This reduces Proposition 4.4.1 to the following construction.

\tm{4.5.1. Proposition} For any integers $k\in\N$ and $g\ge
r\ge1$, there exists a relatively minimal real algebraic Lefschetz
fibration, $p\: W\to\Cp1$, $\conj_W\:W\to W$, with $H_1(W)=0$, and
a real fiber $\SS=p^{-1}(b)$ of genus $g$ such that \roster \item
$\R\SS$ has $r$ components, and $\SS\smallsetminus\R\SS$ is
connected, \item there is a connected orientable component of $\R
W\smallsetminus\R\SS$, whose closure, $F_W$, is bounded by
$\R\SS$, \item $g(F_W)\ge k$.
\endroster
\endtm

\subheading{4.6. Proof of Proposition 4.5.1} Consider a double
covering $q\:W\to\Cp1\times\Cp1$, branched along a non-singular
curve $\C A$ defined over $\R$ and having degree $(2g+2,2d)$,
where $d$ is sufficiently large. The Lefschetz fibration that we
need is the composition of $q$ with the projection to the first
factor. A generic fiber, $\SS_t=q^{-1}(t\times\Cp1)$, $t\in\Cp1$,
projects to $\Cp1$ as a double cover branched at $(2g+2)$ points,
$\C A_t=\C A\cap (t\times\Cp1)$, and thus has genus $g$. If this
branching locus has  $2r$ real points, $\R A_t=\C A_t\cap\Rp1$,
then $\SS_t$ satisfies the condition (1) of Proposition 4.5.1.
 \midinsert
\epsfbox{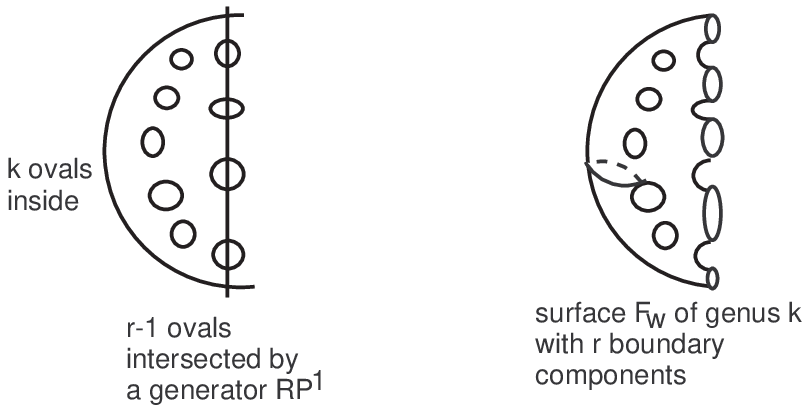} \botcaption{Figure 1}
\endcaption
\endinsert
The conditions (2)--(3) will be obviously satisfied if we choose
the curve $\R A$ and a generator $t\times\Rp1$ with the mutual
position shown on the Figure 1.
The curve with such position can be easily constructed by a
perturbation of a nodal curve in $\Cp1\times\Cp1$ splitting into a
union of generators. \qed

\heading \S5. Some generalizations \endheading

\subheading{5.1. Variants of $S_{X,\SS}$}
 One can consider a version of the invariant $S_{X,\SS}$,
$$S_{X,\SS,K}\:\SC(X,\SS)/K\to R_{X,\SS,K}$$
 which depends on a subgroup $K\subset H^1(\SS)$.
Here $\SC(X,\SS)/K$ is the quotient of $\SC(X,\SS)$ by the action
of $K$, where $h\in K$ acts as $\ss\mapsto \ss+\delta(h)$, and
$\delta\:H^1(\SS)\to H^2(X,\SS)$ is the boundary homomorphism. The
ring $R_{X,\SS,K}$ is
$\A^*_{X,\SS,K}=\Wedge(H^1(X,\SS)/K)\otimes\Z[U]$ in the case of
$S_X=\SW_X$ and $\A^*_{X,\SS,K}\otimes\Z/2$ in the case of
$S_X=\OS_X$. Considering the invariant in the reduced form, we let
$R_{X,\SS,K}$ be just $\Z$ for SW and $\Z/2$ for OS invariants.

The definition of $S_{X,\SS,K}$ is similar to that of $S_{X,\SS}$,
except that the condition $H_1(W)=0$ for a Lefschetz fibration
$W\to S^2$ is replaced by the condition $K=\Im(H^1(W)\to
H^1(\SS))$ (if $K$ can be expressed as such an image). In
particular, for $K=0$, we have $S_{X,\SS,K}=S_{X,\SS}$ and for
$K=H^1(\SS)$ the invariant $S_{X,\SS,K}$ coincides with the
restriction, $S_{X,\ss_\SS}$, of the absolute invariant.

In general, there is a splitting formula
$S_{X,\SS,K}[\ss]=\sum_{\ss'\in[\ss]}S_{X,\SS}(\ss')$ or
equivalently,
$\underline{S}_{X,\SS,K}=(\abs_K)_*(\underline{S}_{X,\SS})$, where
$(\abs_K)_*$ is the push-forward morphism of the projection
$\abs_K\:\SC(X,\SS)\to \SC(X,\SS)/K$. The proof of this formula is
analogous to the proof of (5) in Theorem 1.4.1.

\subheading{5.2. A refinement of the Ozsv\'ath-Szab\'o
$4$-dimensional invariant with respect to a mapping torus boundary
component} The idea used in the definition of $S_{X,\SS}$ can be
used also to define the refinement of the $4$-dimensional
 Ozsv\'ath-Szab\'o invariants in a more general setting. Assume
for instance that $X\: M_0\to M_1$ is a cobordism between
$3$-manifolds, where  $M_1=M_f$ is a mapping torus of some
homeomorphism $f\:\SS\to\SS$.

The plane field tangent to the fibers of the projection $M_f\to
S^1$ defines a canonical $\SC$ extension of the tangent bundle
$\tau_M$. Let $\SC(X,M_f)$ denote the set of the relative $\SC$
structures in $X$ with respect to such $\SC$ extension over
$M_f\subset\partial X$.

Choose any $\rr\in\SC(X,M_1)$ and let $\ss=\abs(\rr)\in\SC(X)$,
and $\ttt_i=\ss|_{M_i}$, $i=0,1$ (here $\ttt_1$ is the canonical
structure on $M_f$). Consider an auxiliary cobordism $W\: M_1\to
M_2$ which has structure of a Lefschetz fibration $q\:W\to
S^1\times[1,2]$ over the annulus, so that $M_i=q^{-1}(S^1\times
i)$. There is a canonical relative $\SC$ structure,
$\rr_{W}\in\SC(W,M_1)$, which is a refinement of the canonical
absolute structure $\ss_W=\abs(\rr_W)\in\SC(W)$. We assume that
the Lefschetz fibration is relatively minimal and $b_1(W)=0$ (one
can always find such a fibration bounded by any prescribed mapping
tori $M_i$, $i=1,2$; for example we may assume in addition that
$M_2\cong\SS\times S^1$).

The homomorphism
$$F^+_{X,\rr}\:HF^+(M_0,\ttt_0)\to
HF^+(M_1,\ttt_1)\cong\Z
$$
is defined as the composition of $F^+_{X\cup W,\rr\rr_W}\:
HF^+(M_0,\ttt_0)\to HF^+(M_2,\ttt_2)$ and the inverse
$(F^+_{W,\ss_W})^{-1}$ to the isomorphism $F^+_{W,\ss_W}\:
HF^+(M_1,\ttt_1)\to HF^+(M_2,\ttt_2)$.

The action of $H_1(X\cup W)=H_1(X\cup W,W)=H_1(X,M_1)$ in
$F^+_{X\cup W,\rr\rr_W}$ composed with the isomorphism
$(F^+_{W,\ss_W})^{-1}$ defines an action  of $\A_{X,M_1}$ in
$F_{X,\rr}$.

The product formula of \cite{OS3} applied to the cobordism $X\cup
W\: M_0\to M_2$ implies a splitting
$F^+_{X,\ss}=\sum_{\rr\in\abs^{-1}(\ss)}F^+_{X,\rr}$. Similarly
one can define maps $F^-_{X,\rr}$ and $F^\mix_{X,\rr}$ and obtain
analogous decompositions of $F^-_{X,\ss}$ and $F^\mix_{X,\ss}$.

All the constructions in this section admit also similar versions
for SW invariants.


\widestnumber \key{ABC} \Refs

\ref
 \key{CW}
 \by A. L. Carey, B.-L. Wang
 \paper Seiberg-Witten-Floer Homology and Gluing Formulae
 \jour Acta Math. Sin.
 \vol19
 \issue2
 \yr2003
 \pages245--296
\endref


\ref
 \key{FS1}
 \by R. Fintushel, R. Stern
 \paper Knots, Links, and $4$-manifolds
 \jour Invent. Math.
 \vol134
 \issue2
 \yr1998
 \pages363--400
\endref

\ref
 \key{FS2}
 \by R. Fintushel, R. Stern
 \paper Surfaces in $4$-manifolds
 \jour Math. Res. Lett.
 \vol4
 \issue
 \yr1997
 \pages907--914
\endref


\ref
 \key{G}
 \by R. Gompf
 \paper Sums of elliptic surfaces
 \jour J. Diff. Geom.
 \vol34
 \issue
 \yr1991
 \pages93--114
\endref


\ref
 \key{KM}
 \by P.B. Kronheimer, T.S. Mrowka
 \paper The genus of embedded surfaces in the projective plane
 \jour Math. Research Letters
 \vol 1
 \issue
 \yr1994
 \pages797--808
\endref

\ref
 \key{MST}
 \by J.W. Morgan, Z.Szab\'o, C.H. Taubes
 \paper A product formula for Seiberg-Witten invariants and the
 generalized Thom Conjecture
 \jour J. Diff. Geom.
 \vol44
 \issue
 \yr1996
 \pages706--788
\endref

\ref
 \key{MW}
 \by V. Munoz, B.-L. Wang
 \paper
 \jour Seiberg-Witten-Floer homology of a surface times a circle
 for non-torsion $\SC$ structures,
ArXive: math.DG/9905050
\endref


\ref
 \key{OS1}
 \by P.~Ozsv\'ath, Z.~Szab\'o
 \paper The symplectic Thom Conjecture
 \jour Ann. of Math.
 \vol151
 \issue1
 \yr2000
 \pages93--124
\endref

\ref
 \key{OS2}
 \by P.~Ozsv\'ath, Z.~Szab\'o
 \paper Higher type adjunction inequalities in Seiberg-Witten
 theory
 \jour J. Differential Geom.
 \vol55
 \issue3
 \yr2000
 \pages385--440
\endref

\ref
 \key{OS3}
 \by P.~Ozsv\'ath, Z.~Szab\'o
 \paper Holomorphic Triangles and invariants for smooth
 four-manifolds
 \jour math.SG/0110169
 \vol
 \issue
 \yr2001
 \pages
\endref

\ref
 \key{OS4}
 \by P.~Ozsv\'ath, Z.~Szab\'o
 \paper Holomorphic triangles invariant and the topology of
 symplectic four-manifolds
 \jour Duke Math. J.
 \vol121
 \issue
 \yr2004
 \pages1--34
\endref


\ref
 \key{T1}
 \by C. Taubes
 \paper The Seiberg-Witten invariants and symplectic forms
 \jour Math. Research. Letters
 \vol1
 \issue
 \yr1994
 \pages809--822
\endref

\ref
 \key{T2}
 \by C. Taubes
 \paper More constraints on symplectic forms from Seiberg-Witten
 invariants
 \jour Math. Research. Letters
 \vol2
 \issue
 \yr1995
 \pages9--13
\endref

\ref
 \key{T3}
 \by C. Taubes
 \paper GR=SW: Counting curves and connections
 \jour J. Diff. Geom.
 \vol52
 \issue
 \yr1999
 \pages453--609
\endref

\ref
 \key{T4}
 \by C. Taubes
 \paper $\roman{SW}\Rightarrow Gr$: from the Seiberg-Witten equations to pseudo-holomorphic curves
\jour J. Amer. Math. Soc.
\vol9
 \issue2
 \pages845--918
\endref

\ref
 \key{T5}
 \by C. Taubes
 \paper The Seiberg-Witten invariants and $4$-manifolds with
 essential tori
 \jour Geometry and Topology
 \vol5
 \issue
 \yr2001
 \pages441--519
\endref


\endRefs
\enddocument